\documentclass[10pt, twoside, a4paper]{article}
\let\finishall\relax\let\Finishall\relax\let\getprepared\relax
\ifx\TestIngCommand\undefined\relax
\else\input ../../../proc/ppreamb-book.tex 
  \input init.tex \fi
\let\TestIngCommand\undefined

\newtheorem{remark}{Remark}

\newtheorem{conj}{Conjecture}
\usepackage{amsmath,amscd,amssymb}                                         
\usepackage{graphics}                                                     
\newtheorem{theo}{Theorem}                                                 
\newtheorem{lem}{Lemma}                                                    
\newtheorem{defi}{Definition}
\newskip\ttglue\ttglue=.5em plus.25em minus.15em                           
\def\firstname#1{\def\FIRSTNAME{#1}\ignorespaces}
\def\lastname#1{\def\LASTNAME{#1}\ignorespaces}
\def\middleinitial#1{\def\MIDDLEINI{#1}\ignorespaces}
\def\department#1{\def\DEPARTMENT{#1}\ignorespaces}
\def\institute#1{\def\INSTITUTE{#1}\ignorespaces}
\def\address#1{\def\ADDRESS{#1}\ignorespaces}
\def\country#1{\def\COUNTRY{#1}\ignorespaces}
\def\otheraffiliation#1{\def\OTHERAFFILIATION{#1}\ignorespaces}
\def\email#1{\def\EMAIL{#1}\ignorespaces}
\newcount\autcount\autcount=0                                              
\newcount\affcount\affcount=0                                              
\newcount\numcount\newcount\nummcount                                      
\newcount\nummmcount\newcount\nummmmcount                                  
\def\writename#1#2{\ \kern-1ex\hbox{
  \csname AUthor\the#1\endcsname\                                          
  \edef\TESTSTR{}\expandafter\ifx\csname auTHor\the#1\endcsname\TESTSTR    
  \else\csname auTHor\the#1\endcsname.\ \fi                                
  \csname authOR\the#1\endcsname$^{\csname AFF\the#1\endcsname}$
  \expandafter\ifx\csname corr\number#1\endcsname\relax                    
  \else\thanks{Corresponding author.}\ \fi                                 
  }\ifnum#1<#2, \else\ \kern-1ex\fi}
\def\writeemail#1{
  \nummcount=0\relax\nummmcount=0\relax                                    
  \loop\ifnum\nummcount<\autcount\advance\nummcount by1\relax              
    {\expandafter\ifnum\csname AFF\the\nummcount\endcsname=#1\relax        
    \global\advance\nummmcount by1\fi}\repeat                              
  \nummcount=0\relax\nummmmcount=0\relax                                   
  \loop\ifnum\nummcount<\autcount\advance\nummcount by1\relax              
    {\expandafter\ifnum\csname AFF\the\nummcount\endcsname=#1\relax        
    \global\advance\nummmmcount by1\relax\def\blank{}\expandafter          
    \ifx\csname EMAIL\the\nummcount\endcsname\blank(no e-mail)
    \else\csname EMAIL\the\nummcount\endcsname                             
    \fi                                                                    
    \ifnum\nummmmcount<\nummmcount; \fi\fi}\repeat}
\long\def\BeginAuthorList#1\EndAuthorList{#1\relax                         
  \author{\vbox{\hsize=390pt\noindent\numcount=0\relax                     
    \loop\ifnum\numcount<\autcount\advance\numcount by1\relax              
      \writename{\numcount}{\autcount}
      \repeat}\\[2mm]                                                      
    \vbox{\small\numcount=0\relax                                          
      \loop\ifnum\numcount<\affcount\advance\numcount by1\relax            
        \vbox{{\count0=\numcount\relax                                     
          \loop\expandafter\ifnum\csname AFF\the\count0\endcsname
            <\numcount\relax\advance\count0 by1\relax\repeat               
          $^{\csname AFF\the\count0\endcsname}$}
        \def\BLANK{}\expandafter\ifx\csname DEPT\the\numcount\endcsname    
          \BLANK                                                           
          \else\csname DEPT\the\numcount\endcsname, \fi                    
        \csname INST\the\numcount\endcsname,                               
        \csname ADDR\the\numcount\endcsname,                               
        \csname COUN\the\numcount\endcsname                                
        \edef\TEST{}\expandafter\ifx\csname OTHE\the\numcount\endcsname
          \TEST                                                            
          .\else;\break\csname OTHE\the\numcount\endcsname.\fi}
        \vbox{\writeemail{\numcount}}
        \repeat}\\}}
\expandafter\def\csname x1\endcsname{}
\expandafter\def\csname x2\endcsname{}
\expandafter\def\csname x3\endcsname{}
\expandafter\def\csname x4\endcsname{}
\expandafter\def\csname x5\endcsname{}
\expandafter\def\csname x6\endcsname{}
\expandafter\def\csname x7\endcsname{}
\expandafter\def\csname x8\endcsname{}
\expandafter\def\csname x9\endcsname{}
\def\Author#1#2{\global\advance\autcount by1\relax#2                       
  \expandafter\edef\csname AUthor\the\autcount\endcsname{\FIRSTNAME}
  \expandafter\edef\csname auTHor\the\autcount\endcsname{\MIDDLEINI}
  \expandafter\edef\csname authOR\the\autcount\endcsname{\LASTNAME}
  \expandafter\edef\csname EMAIL\the\autcount\endcsname{\EMAIL}
  \let\tempera\"\def\"{\string\"}\expandafter\ifx\csname x\DEPARTMENT
    \endcsname\relax                                                       
    \global\advance\affcount by1\relax\let\"\tempera                       
    \expandafter\edef\csname DEPT\the\affcount\endcsname{\DEPARTMENT}
    \expandafter\edef\csname INST\the\affcount\endcsname{\INSTITUTE}
    \expandafter\edef\csname ADDR\the\affcount\endcsname{\ADDRESS}
    \expandafter\edef\csname COUN\the\affcount\endcsname{\COUNTRY}
    \expandafter\edef\csname OTHE\the\affcount\endcsname{\OTHERAFFILIATION}
    \expandafter\edef\csname AFF\the\autcount\endcsname{\the\affcount}
  \else\expandafter\edef\csname AFF\the\autcount\endcsname{\DEPARTMENT}
  \fi\let\"\tempera\ignorespaces}
\def\CorrespondingAuthor#1#2{
  \expandafter\xdef\csname corr\number#1\endcsname{cor}
  \Author#1{#2}}
\def\PaperTitle#1{\title{\bf#1}}
\def\Category#1{\ignorespaces}
\def\keywords#1{{\noindent \emph{Keywords:}                                
  \def\BLANK{}\def\TEST{#1}\ifx\BLANK\TEST(n/a).\else#1\fi}}
\getprepared                                                               
\begin{document}                                                           

\PaperTitle{Filters and Functions in Multi-scale Constructions: Extended Abstract}%
\Category{(Pure) Mathematics}

\date{}

\BeginAuthorList
  \Author1{
    \firstname{Wayne}
    \lastname{Lawton}
    \middleinitial{M}   
    \department{Applied Mathematics Program, Science Division}
    \institute{Mahidol University International College}
    \address{Nakornpathom 73170}
    \country{Thailand}
    \otheraffiliation{School of Mathematics and Statistics, University of Western Australia, Perth, Australia}
    \email{wayne.lawton@uwa.edu.au}}
\EndAuthorList
\maketitle
\thispagestyle{empty}
\begin{abstract}
We derive results about geometric means of the Fourier modulus of filters and functions
related to refinable distributions with arbitrary dilations and translations. Then we
develop multi-scale constructions for dilations by Pisot-Vijayaraghavan numbers and 
translations in associated quasilattices.
\end{abstract}

\finishall
\section{Introduction and Statement of Results}
\label{intro}
In this paper $\mathbb{Z}, \mathbb{N} = \{1,2,3,...\}, \mathbb{Q}, \mathbb{R}, \mathbb{C}$ denote the integer, natural, rational, real, and complex numbers and $\mathbb{T} = \mathbb{R}/\mathbb{Z}$ denotes the circle group.
\\ \\
Let $\lambda \in \mathbb{R}\, \backslash \, [-1,1],$ $a_1,....a_m \in \mathbb{C}\backslash\{0\},$
$a_1+\cdots+a_m = |\lambda|,$ $\tau_1 < \cdots < \tau_m \in \mathbb{R},$ and
\begin{equation}
\label{A}
A(y) = |\lambda|^{-1} \sum_{j=1}^{m} a_j\, e^{2\pi i \tau_j y}, \ \ y \in \mathbb{R}.
\end{equation}
There exists a unique compactly supported distribution $f$ that
satisfies $\int f(x)dx = 1$ and the refinement equation
\begin{equation}
\label{refinable}
    f(x) = \sum_{j=1}^m a_j f(\lambda \, x - \tau_j)
\end{equation}
The Fourier transform ${\widehat f}(y) = \int_{-\infty}^{\infty} \, f(x)e^{2\pi i x y} \, dx$ satisfies ${\widehat f}\,(0) = 1$ and admits the product expansion
\begin{equation}
\label{product}
    {\widehat f}\,(y) = \prod_{k=1}^{\infty} A\left(\, y\lambda^{-k}\, \right).
\end{equation}
$A(y)$ and ${\widehat f}\,(y)$ are the filters and functions of this paper. Equation \ref{product} implies that they are related by the functional equation
\begin{equation}
\label{functional}
    {\widehat f}\, \left(\, y \, \lambda^k \, \right) = {\widehat f}\,(y) \,
    \prod_{j = 1}^{k} A\left(\, y \, \lambda^j\, \right), \ \ k \in \mathbb{N}.
\end{equation}
The Mahler measure \cite{mahler} (or height) $M(P)$ of a polynomial $P(z_1,...,z_d)$ is defined by
\begin{equation}
\label{mahler}
    M(P) = \exp \int_{t_1=0}^{1} \cdots \int_{t_d=0}^{1}
            \ln |P(e^{2\pi i \, t_1},...,e^{2\pi i \, t_d})|\, dt_1 \dots dt_d.
\end{equation}
For $P(z) = c(z-\omega_1)\cdots(z-\omega_q)$ Jensen's theorem \cite{jensen} gives
$M(P) = |c| \prod_{j = 1}^{q} \max\,  \{\, 1,|\,\omega_j\,|\, \}.$
Mahler measure arises in prediction theory of stationary
random processes \cite{lawton1}, algebraic dynamics \cite{everestward}, and
in Lehmer's open problem in number theory \cite{lehmer}.
\\ \\
$A(y)$ is an almost periodic function in the (uniform) sense of Bohr \cite{bohr}.
Let $d$ be the maximal number of linearly independent $\tau_1,...,\tau_m$ over $\mathbb{Q}.$
$A(y)$ is periodic if and only if $d = 1.$
There exists a polynomial $P(z_1,...,z_d)$ and real numbers $r_1,...,r_d$ such that
\begin{equation}
\label{AP}
    A(y) = P(e^{2\pi i r_1y},...,e^{2\pi i r_dy}), \ \ y \in \mathbb{R}.
\end{equation}
\begin{theo}
\label{theorem1}
$\ln |A(y)|$ is almost periodic in the $($mean squared$\,)$ sense of Besicovitch \cite{besicovitch}.
Furthermore
\begin{equation}
\label{mean}
\lim_{L \rightarrow \infty} \frac{1}{2L} \int_{-L}^{L} \ln |A(y)| = \ln M(P).
\end{equation}
\end{theo}
\begin{defi}
\label{MA}
The Mahler measure of $A$ is $M(A) = M(P).$
\end{defi}
\begin{theo}
\label{theorem2}
    Define $\rho(f) = - \ln M(A)/\ln \lambda.$
\begin{equation}
\label{Asympt1}
    \lim_{L \rightarrow \infty} \frac{1}{2L \ln L} \int_{-L}^{L} \ln |{\widehat f}\, (y)|\, dy = - \rho(f).
\end{equation}
\end{theo}
The Hausdorff dimension of a Borel measure is defined in (\cite{peresschlagsolomyak}, p. 6) as
\begin{equation}
\label{HD}
    \hbox{dim}(\nu) = \inf \{ \, \hbox{dim}(B) \, : \, B \hbox{ is a Borel subset of } \mathbb{R} \hbox{ and }
    \nu (\mathbb{R} \backslash B) = 0 \ \}
\end{equation}
where $\hbox{dim}(B)$ is the Hausdorff dimension of $B.$
\begin{conj}
\label{hausdim}
If $f$ is a Borel measure then $\rho(f) = \hbox{dim}(f).$
\end{conj}
If conjecture \ref{hausdim} is validated we propose to call $\rho(f)$ the Hausdorff dimension of $f.$
\\ \\
An algebraic integer $\lambda$ is a Pisot–-Vijayaraghavan (PV) number if the roots $\lambda = \lambda_1,...,\lambda_n$ of its minimal polynomial $\Lambda(z)$ satisfy
$
    |\lambda_j| < 1, \, j = 2,...,n.
$
They were discovered by Thue \cite{thue}, rediscovered by Hardy \cite{hardy} who with Vijayaraghavan studied their Diophantine approximation \cite{cassels} properties, and studied
by Pisot \cite{pisot} in his dissertation. Examples are:
$$
\begin{array}{ccccc}
\lambda =\lambda_1  & \lambda_2 & \lambda_2 & \lambda_4 & \hbox{minimal polynomial} \\
1.6180 &-0.6180 &  &  &  z^2 - z - 1 \\
-1.6180 & 0.6180 & & & z^2 + z - 1 \\
3.4142 &  0.5858 &  &  &  z^2 - 4z + 2 \\
2.2470 & 0.5550 & -0.8019 & & z^3 - 2z^2 -z + 1 \\
1.3247 & -0.6624 + 0.5623i & -0.6624 - 0.5623i &  &  z^3-z-1\\
1.3803 &  -0.8192 & 0.2194 + 0.9145i  &  0.2194 - 0.9145i &  z^4-z^3 -1
\end{array}
$$
If $\lambda$ is a PV number, $m = 2,$ $a_1 = a_2 = |\lambda|/2,$ $\tau_1 = 0,$ and $\tau_2 = 1$ then $f$ is a measure defined by a Bernoulli convolution and Erd\"{o}s \cite{erdos} proved that there exists $\alpha > 0$ and $\gamma \in \mathbb{C} \backslash \{0\}$ such that
\begin{equation}
\label{erdos}
    \lim_{k \rightarrow \infty} {\widehat f}\, (\alpha \lambda^k) = \gamma.
\end{equation}
Erd\"{o}s' proof used Equation \ref{functional} and the fact that the
set of real $\alpha,$ for which there exists a sequence of integers $n_k$ with $|\alpha \, \lambda^k - n_k| \rightarrow 0$ exponentially fast, is dense.
\\ \\
An algebraic integer $\lambda$ is a Salem number \cite{salem1} if the roots $\lambda = \lambda_1,...,\lambda_n$ of its minimal polynomial satisfy
$
    |\lambda_j| \leq 1, \, j = 2,...,n.
$
Like PV numbers, Salem numbers appear in Diophantine approximation and harmonic analysis \cite{meyer,salem2}. The smallest $1.17628$ is the Mahler measure of Lehmer's polynomial
$z^{10}+z^9-z^7-z^6-z^5-z^4-z^3+z+1$
and is conjectured to be the smallest known Mahler measure $> 1$ of any polynomial.
\\ \\
Kahane \cite{kahane} extended Erd\"{o}s' result for $\lambda$ a Salem number and Dai, Feng and Wang \cite{daifengwang} extended Kahane's result for $m > 2$ and $\tau_j \in \mathbb{Z}.$
\begin{theo}
\label{theorem3}
If $\lambda$ is a PV number with degree $n,$ $\alpha \in \mathbb{Q}(\lambda),$ $\{\tau_1,...,\tau_m\} \subset \mathbb{Z}\left[\, 1,\lambda,...,\lambda^{n-1}\, \right],$ and ${\widehat f}\, (\alpha \, \lambda^k) \neq 0$ for all $k \in \mathbb{N}$ then the following limit exists
\begin{equation}
\label{Asympt2}
    \lim_{k \rightarrow \infty} \frac{\ln|{\widehat f}\, (\alpha \, \lambda^k)|}{k\, \ln |\lambda|}
\end{equation}
and equals the mean value of $\ln |P_{\hbox{trig}}|$ over a finite orbit of $C$ in the $n$-dimensional torus group. Here $P_{\hbox{trig}}$ is the trigonometric polynomial on the torus constructed
from the polynomial $P$ by Equation \ref{Ptrig} and the polynomial $P$ is related to $A(y)$ by
Equation \ref{AP}. If $\lambda$ is not an integer then $\alpha$ can be chosen such that $\lim_{k\rightarrow}{\widehat f}\, (\alpha \lambda^k) = \gamma$ where $\gamma \in \mathbb{C}\backslash \{0\}$ thus extending the results in \cite{daifengwang} to noninteger values of $\tau_j.$
\end{theo}
Section \ref{prelim} reviews notation and analytical results
for Theorems \ref{theorem1} and \ref{theorem2}, which are proved in Section \ref{proofs},
and algebraic results for Theorem \ref{theorem3}, which is proved in Section \ref{proofs},
and for Lemma \ref{QL3}, which is proved in Section \ref{multiscale}.
Section \ref{multiscale} constructs multi-scale analyses for spaces of distributions
based on translations by points in quasilattice subsets of $\mathbb{R}.$
\section{Preliminary Results}
\label{prelim}
Let $\mu$ be Lebesgue measure. For a function $K : I \rightarrow \mathbb{C}$ on an interval $I \subseteq \mathbb{R}$ let $K^{\prime}$ be its derivative.
For $n \in \mathbb{N},$ $\mathbb{R}^n$ is $n$-dimensional Euclidean vector space with norm $|| \, ||,$ $\mathbb{Z}^n$ is a rank $n$ lattice subgroup, $\mathbb{T}^n =  \mathbb{R}^n/\mathbb{Z}^n$ is the $n$-dimensional torus group equipped with Haar measure, and $p_n : \mathbb{R}^n \rightarrow \mathbb{T}^n$ is the canonical epimorphism. Vectors are represented as column vectors and $T$ denotes transpose. We define $|| \ ||_{\, \mathbb{T}^n}: \mathbb{T}^n \rightarrow \mathbb{R}$ by
$$||\,g\,||_{\, \mathbb{T}^n} = \min \{ \, ||\,x\,|| : x \in \mathbb{R}^n \hbox{ and } p_n(x)=g \, \}, \ \ g \in \mathbb{T}^n.$$
\begin{lem}
\label{titchmarsh}
$(${\bf Titchmarsh}$)$
If $h$ is a compactly supported distribution and $[a,b]$ is
the smallest interval containing the support of $h$ then
\begin{equation}
\label{exptype}
    H(z) = \int_{a}^{b} h(t) e^{2\pi zt} dt
\end{equation}
is an entire function of exponential type. Moreover the number $n(r)$
of zeros of $H(z)$ in the disk $|z| < r$ satisfies
\begin{equation}
\label{titchmarsheqn}
    \lim_{r \rightarrow \infty} \frac{n(r)}{r} = 2(b-a).
\end{equation}
\end{lem}
Titchmarsh's proof (Theorem IV, \cite{titchmarsh}) assumed that $h$ was a Lebesgue integrable function but a simple regularization argument implies that if $h$ is a compactly supported distribution then
Equation \ref{titchmarsheqn} holds. Therefore there exists $\beta > 0$ such that
\begin{equation}
\label{nzeros}
    \hbox{card } \{ \, y \, : \, y \in [-L,L] \hbox{ and } H(y) = 0 \, \}  \leq \beta\, L, \ \ L \geq 1
\end{equation}
\begin{lem}
\label{uvbound}
Let $I = [a,b]$ and $K : I \rightarrow \mathbb{C}$ be differentiable
such that the image $K^{\prime}(I)$ of $I$ under $K^{\prime}$ is contained
in a single quadrant of $\mathbb{C}.$ If $u >0, v > 0$ and $|K| \leq v$
over $I$ and $|K^{\prime}| \geq u$ over $I$ then
\begin{equation}
\label{muI}
    \mu(I) \leq \frac{2{\sqrt 2}v}{u}.
\end{equation}
\end{lem}
{\bf Proof} We present the proof that we gave in (\cite{lawton3}, Lemma 1).
Since $u \leq |K^{\prime}|$ over $I$ and $\mu(I) = b-a,$ the triangle
inequality
$|K^{\prime}| \leq |\Re\, K^{\prime}| + |\Im\, K^{\prime}|$
gives
$$u \, \mu(I) \leq \int_{a}^{b} |K^{\prime}(y)| \, dy \leq \int_{a}^{b}
\left( \, |\Re \, K^{\prime}(y)| +  |\Im \, K^{\prime}(y)| \, \right) \, dy.$$
Since $K^{\prime}(I)$ is contained in a single quadrant of $\mathbb{C}$ there exist
$c \in \{1,-1\}, d \in \{1,-1\}$ such that $|\Re \, K^{\prime}(y)| = c\, \Re \, K^{\prime}(y)$
and $|\Im \, K^{\prime}(y)| = d \, \Im \, K^{\prime}(y)$ for all $y \in I.$
Therefore
$$\int_{a}^{b}
\left( \, |\Re \, K^{\prime}(y)| +  |\Im \, K^{\prime}(y)| \, \right) \, dy
= (c\, \Re K(b) + d\, \Im K(b)) - (c\, \Re K(a) + d\, \Im  K(a)).$$
Since $|c| = |d| = 1$ and $|K(y)| \leq v$ over $I,$ the quantity on right
is bounded above by $2{\sqrt 2} \, v.$ Combining these three inequalities completes the proof.
\begin{lem}
\label{boyd}
There exists a sequence $c_2,c_3,...> 0$ such that for all monic polynomials $P(z)$ with $k \geq 2$ nonzero coefficients
\begin{equation}
\label{pb}
    \mu(\{\, y \in [0,1] \, : \, |P(e^{2\pi i \, y})| \leq v \, \}) \leq c_k\, v^{1/(k-1)}, \ \ v > 0.
\end{equation}
\end{lem}
{\bf Proof} We proved this in (Theorem 1, \cite{lawton3}) using induction on $k,$ lemma \ref{uvbound}, and the fact that a polynomial of degree $q$ can have at most $q$ distinct roots.
\\ \\
{\bf Remark} We used Lemma \ref{boyd} and the method we developed in \cite{lawton2} and Boyd
developed in \cite{boyd} to prove a conjecture that he formulated in \cite{boyd}. The
conjecture expresses the Mahler measure of a multidimensional polynomial as a limit of
Mahler measures of univariate polynomials and provides an alternative proof of the
special case of Lehmer's conjecture proved in \cite{dobrowolski}. The conjecture is stated
in http:$//$en.wikipedia.org$/$wiki$/$Mahler$\textunderscore$measure
\begin{lem}
\label{tb}
There exists a sequence $c_2,c_3,...> 0$ such that
\begin{equation}
\label{tb1}
    \mu(\{\, y \in [-L,L] \, : \, |A(y)| \leq v \, \})
    \leq L\, c_m \, v^{1/(m-1)}, \ \ L \geq 1, \ v > 0.
\end{equation}
\end{lem}
{\bf Proof} We observe that Lemma \ref{tb} concerns (not necessarily periodic) trigonometric polynomials. The conclusion is obvious for $m = 2$ so we assume that $m \geq 3$ and proceed by induction on $m$ using a technique similar to the one we used to prove Lemma \ref{boyd}. Since $|e^{-2\pi i y}A(y)| = |A(y)|$ we may assume without loss of generality that $0 = \tau_1 < \cdots \tau_m.$ Then $A^{\prime}(y)$ is a trigonometric polynomial with $m-1$ terms
so by induction there exists $c_{m-1} > 0$ such that
\begin{equation}
\label{tb2}
    \mu(\{\, y \in [-L,L] \, : \, |A^{\prime}(y)| \leq u \, \})
    \leq L\, c_{m-1} \, u^{1/(m-2)}, \ \ L \geq 1, \, u > 0.
\end{equation}
We observe that $A$ and $A^{\prime},$ their real and imaginary parts, and their squared moduli are restrictions of functions of exponential type so they satisfy the hypothesis in Lemma \ref{titchmarsh}. Therefore Equation \ref{nzeros} implies that there exists $\beta > 0$ such that for all $u >0, v > 0$ and $L \geq 1$ the subset of $[-L,L]$ where $|A^{\prime}(y)| \geq u$ and $|A(y)| \leq v$ can be expressed as the union of not more than $\beta\, L$ closed intervals $I$ such that $A^{\prime}(I)$ is contained in a single quadrant of $\mathbb{C}.$ Therefore Lemma \ref{uvbound} implies that there exists $\beta > 0$ such that
\begin{equation}
\label{tb3}
    \mu(\{\, y \in [-L,L] \, : \, |A(y)| \leq v \hbox{ and } |A^{\prime}(y)| \geq u \, \})
    \leq \frac{L\, \beta \, v}{u}, \ \ L \geq 1, \, u > 0, \, v > 0.
\end{equation}
Combining Equations \ref{tb2} and \ref{tb3} gives for every $L \geq 1,\, v > 0$
\begin{equation}
\label{tb4}
    \mu(\{\, y \in [-L,L] \, : \, |A(y)| \leq v  \, \}) \leq L \min_{u > 0}
        \left[\, c_{m-1} \, u^{1/(m-2)} + \frac{\beta \, v}{u} \, \right] = Lc_m v^{1/(m-1)}
\end{equation}
where
$c_{m} = (1+\beta(m-2))\left[\, \frac{c_{m-1}}{\beta(m-2)}\, \right]^{\frac{m-2}{m-1}}.$
\\ \\
For $d \in \mathbb{N}$ and $r = [r_1,...,r_d]^T \in \mathbb{R}^d$ we define the homomorphism
$\Psi_r \, : \, \mathbb{R} \rightarrow \mathbb{T}^d$
by
\begin{equation}
\label{psi}
    \Psi_r(y) = p_d(y r), \ \ y \in \mathbb{R}
\end{equation}
where $p_d : \mathbb{R}^d \rightarrow \mathbb{T}^d$ is the canonical epimorphism.
\begin{lem}
\label{BSW}
$(${\bf Bohl--Sierpinski--Weyl}$)$ The image of $\Psi_r$ is a dense subgroup of $\mathbb{T}^d$ if and only if the components $r_1,...,r_d$ of $r$ are linearly independent
over $\mathbb{Q}.$ If they are independent then for every continuous function $S : \mathbb{T}^d \rightarrow \mathbb{C}$
\begin{equation}
\label{unifdist}
    \lim_{L \rightarrow \infty} \frac{1}{2L} \, \int_{-L}^{L} S \circ \Psi_r\, (y)dy = \int_{\mathbb{T}^d} S(g)\, dg
\end{equation}
where $S \circ \Psi_r : \mathbb{R} \rightarrow \mathbb{C}$ is the composition of $S$ with $\Psi_r$
and $dg$ is Haar measure on $\mathbb{T}^d.$
\end{lem}
{\bf Proof} Arnold (\cite{arnold}, p. 285--289) discuss the significance and gives Weyl's proof of
this classic result called the {\it Theorem on Averages} and asserts it ``may be found implicitly
in the work of Laplace, Lagrange, and Gauss on celestial mechanics'' and ``A rigorous proof was
given only in 1909 by P. Bohl, V. Sierpinski, and H. Weyl in connection with a problem of
Lagrange on the mean motion of the earth's perihelion."
\\ \\
The analytic results above suffice to prove Theorems \ref{theorem1} and \ref{theorem2}.
The following algebraic results are required to prove Theorem \ref{theorem3} and Lemma \ref{QL3}.
\\ \\
If $\lambda$ is an algebraic integer with minimal polynomial
$\Lambda(z) = z^n + c_{n-1}z^{n-1} + \cdots + c_0, \ c_j \in \mathbb{Z}$
that has roots $\lambda = \lambda_1, \lambda_2, \dots, \lambda_n,$ then the
companion matrix
$$
    C =
    \left[\begin{array}{ccccc}
        0 & 1 & 0 & \hdots & 0 \\
        0 & 0 & 1 & \ddots & \vdots \\
        \vdots & \ddots & \ddots & \ddots & 0 \\
        0 & 0 & \hdots & 0 & 1 \\
        -c_0 & -c_1 & \hdots & \hdots & -c_{n-1}
    \end{array}\right],
$$
Vandermonde matrix
$$
    V =
    \left[\begin{array}{cccc}
     1                 &  1                & \hdots &  1 \\
     \lambda_1         & \lambda_2         & \hdots & \lambda_{n} \\
     \vdots            & \vdots            & \hdots & \vdots \\
     \lambda_{1}^{n-1} & \lambda_{2}^{n-1} & \hdots & \lambda_{n}^{n-1} \\
    \end{array} \right],
$$
and diagonal matrix $D = \hbox{diag}(\lambda_1,...\lambda_n)$ satisfy
\begin{equation}
\label{keyeq}
    C^k \, V = V \, D^k, \ \ k \in \mathbb{N}.
\end{equation}
The matrix $VV^T$ is invertible since $\det VV^T = |c_0|^2$
and it has integer entries since
\begin{equation}
\label{keyeq2}
    (VV^T)_{i,j} = \sum_{\ell = 1}^{n} \lambda_{\ell}^{i+j-2} = \hbox{trace } C^{i+j-2}, \ \ i,j \in \{1,...,n\}.
\end{equation}
Define $v = [\, 1, \, \lambda_1, \, ... \, , \lambda_{1}^{n-1}\, ]^T.$
For $q \in \mathbb{Q}^n$ let
$u = V^T(VV^T)^{-1}q = [\, \alpha_1, \, ... \, , \alpha_n]^T.$
Then $\alpha_j \in \mathbb{Q}(\lambda_j), \, j = 1,...,n$ and $V \, u = q.$ If $\alpha \in \mathbb{Q}$ then there exists $b \in \mathbb{Q}^T$ such that $\alpha = v^T b.$ Then $u = V^Tb = [\alpha_1\, ...\, \alpha_n]^T$ satisfies $\alpha = \alpha_1,$ $\alpha_j \in \mathbb{Q}(\lambda_j),\, j = 1,...,n$ and $q = Vu \in \mathbb{Q}^n.$
\\ \\
Now assume that $\lambda$ is a PV number. Then Equation \ref{keyeq} implies that
\begin{equation}
\label{conv1}
    ||\, C^k q - \alpha \, \lambda^k \, v|| \rightarrow 0
\end{equation}
where the convergence is exponentially fast. Then Equation \ref{conv1}
implies that
\begin{equation}
\label{conv2}
    ||\, p_n(C^k q - \alpha \, \lambda^k \, v)||_{\mathbb{T}^n} = ||\, C^k p_n(q) - p_n(\alpha \, \lambda^k \, v)||_{\mathbb{T}^n} \rightarrow 0.
\end{equation}
Since $p_n(\mathbb{Q}^n)$ is the set of preperiodic points of the endomorphism $C : \mathbb{T}^n \rightarrow \mathbb{T}^n,$ if $\alpha \in \mathbb{Q}(\lambda)$ then the sequence $p_n(\alpha \, \lambda^k \, v) \in \mathbb{T}^n$ converges exponentially fast to a periodic orbit of
$C.$ Conversely every periodic orbit is the limit of such a sequence.
\begin{remark}
We suggest further study of the preperiodic points of $C$ using results from integral matrices discussed by Newman \cite{newman}, the open dynamical version of the Manin-Mumford Conjecture solved by Raynaud \cite{raynaud1}, \cite{raynaud2}, and connections between Mahler
measure and torsion points developed by Le \cite{le}. We also suggest extensions obtained by replacing $\mathbb{R}^n$ by certain Lie groups $($stratified nilpotent groups with rational structure constants$)$ considered in \cite{lawton4}.
\end{remark}
The following result, proved by Minkowsky in 1896 \cite{minkowski}, is the foundation of
the geometry of numbers (\cite{gruberlekkerkerker}, II.7.2, Theorem 1), (\cite{maninpanchishkin},
Theorem 4.7)
\begin{lem}
\label{minkowski}
$(${\bf Minkowski}$)$ If $X \subset \mathbb{R}^n$ is convex, $X = -X,$ and the volume of $X$ exceeds
$2^n$ then $X$ contains a point in $\mathbb{Z}^n\backslash \{0\}.$
\end{lem}
\section{Proofs}
\label{proofs}
For every nonzero polynomial $P(z_1,...\, ,z_d)$
we define the trigonometric polynomial \\
$P_{\hbox{trig}} : \mathbb{T}^d \rightarrow \mathbb{C}$
by
\begin{equation}
\label{Ptrig}
    P_{trig}(p_d([t_1,...\, ,t_d]^T)) = P(e^{2\pi i\, t_1},...\, ,e^{2\pi i\, t_d}).
\end{equation}
and observe through a simple induction based computation that
    $\ln |P_{trig}| \in L^2(\mathbb{T}^d).$
\\
Henceforth we assume that $r = [r_1,...\, ,r_d]^T \in \mathbb{R}^d$ and
$P(z_1,...\, ,z_d)$ are chosen to satisfy Equation \ref{AP}. Therefore
$r_1,...\, ,r_d$ are linearly independent over $\mathbb{Q}$ and moreover
\begin{equation}
\label{AP1}
    A(y) = P_{trig} \circ \Psi_r (y), \ \ y \in \mathbb{R}.
\end{equation}
{\bf Proof of Theorem \ref{theorem1}}
For all $v > 0$ define $S_v : \mathbb{T}^n \rightarrow \mathbb{R}$ by
$$
    S_v(g) = \ln \left( \max \{v,|P_{trig}(g)|\} \right).
$$
Lebesgue's dominated convergence theorem implies that
\begin{equation}
\label{LDC}
    \lim_{v \rightarrow 0} \int_{\mathbb{T}^d} |\, S_v(g)-\ln |P_{trig}(g)| \, |^2 dg = 0.
\end{equation}
Lemma \ref{tb} implies that
\begin{equation}
\label{Av}
    \lim_{v \rightarrow 0} \lim_{L \rightarrow \infty} \frac{1}{2L}
    \int_{-L}^{L} |\, \ln |A(y)| - S_v \circ \Psi_r(y)) \, |^2 dy = 0
\end{equation}
and hence $\ln |A(y)|$ is almost periodic in the sense of Besicovitch.
We combine Equations \ref{LDC} and \ref{Av} with Lemma \ref{BSW} to complete the proof
by computing
$$
\begin{array}{ccccc}
    \lim_{L \rightarrow \infty} \frac{1}{2L} \int_{-L}^{L} \ln |A(y)| dy & = &
    \lim_{v \rightarrow 0} \lim_{L \rightarrow \infty} \frac{1}{2L} \int_{-L}^{L} S_v \circ \Psi_r(y)) dy & = & \\ \\
    \lim_{v \rightarrow 0} \int_{\mathbb{T}^d} S_v(g) dg & = &
    \int_{\mathbb{T}^d} \ln |P_{trig}(g)| dg & = & M(P).
\end{array}
$$
{\bf Remark} Since the zero set of $P_{trig}$ is a real analytic set, an alternative proof
based on Lojasiewicz's structure theorem \cite{kranzparks}, \cite{lojasiewicz} for real analytic sets may be possible using methods that we developed in \cite{lawton5} to prove the Lagarias-Wang Conjecture.
\\ \\
{\bf Proof of Theorem \ref{theorem2}} Let
$L = b|\lambda|^k$ with $k \in \mathbb{N}$ and $b \in [\, |\lambda^{-1}|,1\, )$ then use the functional Equation \ref{functional} and Theorem \ref{theorem1} to compute
$$
\begin{array}{ccc}
\lim_{L \rightarrow \infty} \frac{1}{2L \ln L} \int_{-L}^{L} \ln |{\widehat f}\, (y)|\, dy
& = & \\ \\
\lim_{k \rightarrow \infty} \frac{1}{2 kb\, |\lambda|^k \ln |\lambda|}
\sum_{j = 1}^{k-1} \int_{b|\lambda|^j}^{b|\lambda|^{j+1}} \ln \left(\, |{\widehat f}\, (y){\widehat f}\, (-y)|\, \right)\, dy
& = & \\ \\
\lim_{k \rightarrow \infty} \frac{1}{2 k\, |\lambda|^k \ln |\lambda|}
\sum_{j = 1}^{k-1} |\lambda|^j \int_{1}^{|\lambda|}
\ln \left(\, |{\widehat f}\, (ub|\lambda|^j){\widehat f}\, (-ub|\lambda|^j)|\,\right)\,  du
& = & \\ \\
\lim_{k \rightarrow \infty} \frac{1}{2 k\, |\lambda|^k \ln |\lambda|}
\sum_{j = 1}^{k-1} |\lambda|^j \int_{1}^{|\lambda|}
\ln \left(\, |\, {\widehat f}\, (ub\lambda^j)\, {\widehat f}\, (-ub\lambda^j)\, |\,\right)\, du
& = & \\ \\
\lim_{k \rightarrow \infty} \frac{1}{2 k\, |\lambda|^k \ln |\lambda|}
\sum_{j = 1}^{k-1} |\lambda|^j \int_{1}^{\lambda} \left(\ln |{\widehat f}\, (bu){\widehat f}\, (-bu)| + \sum_{i=1}^j \ln |A(ub\lambda^j)A(-ub\lambda^j)|\right)du
& = & \\ \\
\lim_{k \rightarrow \infty} \frac{1}{2 k\, |\lambda|^k \ln |\lambda|}
\sum_{j = 1}^{k-1} |\lambda|^j \sum_{i=1}^j \int_{1}^{\lambda} \left(
\ln |A(ub\lambda^j)| + \ln |A(-ub\lambda^j)|\right)\, du
& = & \\ \\
\lim_{k \rightarrow \infty} \frac{(|\lambda| - 1)\, \ln M(A)}{k\, |\lambda|^k \, \ln |\lambda|}
\sum_{j = 1}^{k-1} j\, |\lambda|^j
& = & \\ \\
\lim_{k \rightarrow \infty} \frac{(|\lambda| - 1)\, \ln M(A)}{k\, |\lambda|^k \, \ln |\lambda|}
\left(\, -|\lambda| (|\lambda| - 1)^{-2}(|\lambda|^k -1) + k|\lambda|^k(|\lambda|-1)^{-1}\, \right)
& = & \\ \\
\frac{\ln M(A)}{\ln \lambda}.
\end{array}
$$
{\bf Applications of Theorem \ref{theorem2}} The following computations support Conjecture \ref{hausdim}.
\\ \\
The characteristic function $f$ of $[0,1]$ satisfies $f(x) = f(2x) + f(2x-1).$ Therefore
$A(y) = (1 + e^{2\pi i y})/2$ so $\rho(f) = 1$ which equals the Hausdorff dimension of $f.$
We observe that $|{\widehat f}\, (y)| = |\sin y / y|$ decays like $|y|^{-1}.$
\\ \\
The uniform measure $f$ on Cantor's ternary set satisfies $f(x) = \frac{3}{2}f(3x) + \frac{3}{2}f(3x-2).$ Therefore $A(y) = (1 + e^{6\pi i y})/2$ so
$\rho(f) = \ln 2/\ln 3 \approx 0.6309$ which equals the Hausdorff dimension of $f.$
This suggests that ${\widehat f}(y)$ decays like $|y|^{-\ln 2/\ln 3}.$
\\ \\
The measure $f$ defined by a Bernoulli convolution and PV number $\lambda$ satisfies \\
$f(x) = \frac{|\lambda|}{2}f(\lambda x) + \frac{|\lambda|}{2}f(\lambda x - 1).$
Therefore $A(y) = (1 + e^{2\pi i \lambda y})/2$ so $\rho(f) = \ln 2/\ln |\lambda|$ which
equals the Hausdorff dimension of $f$ computed by Peres, Schlag and Solomnyak \cite{peresschlagsolomyak}. This suggests that
${\widehat f}(y)$ decays like $|y|^{-\ln 2/\ln |\lambda|}.$
\\ \\
Let $f$ be a nonzero distribution that satisfies $f(2x) = 2f(2x) + 2f(2x-1) - 2f(2x-3).$
Therefore $A(y) = 1+e^{2\pi i y} - e^{4\pi i y}$ so $\rho(f) = -\ln ((1+\sqrt 5)/2)/\ln 2 \approx 0.6942.$ This suggests that $|{\widehat f}(y)|$ grows like $|y|^{-\rho(f)}$
\\ \\
{\bf Proof of Theorem \ref{theorem3}} Uses algebraic results in Section \ref{prelim} with
the use of the functional Equation \ref{functional} as in the proof of Theorem \ref{theorem2}.
Details will be provided in the full paper.
\section{Multiscale Constructions}
\label{multiscale}
Let $\lambda$ be a PV number and let $C, D, V, \lambda = \lambda_1,...,\lambda_n$ be defined as in
Section \ref{prelim}. A vector $\sigma \in \mathbb{R}^n$ is called admissible if
$\sigma_1 = 0,$ $\sigma_j > 0$ for $j \geq 2,$ and whenever $2 \leq j < k \leq n$
and $\lambda_j$ and $\lambda_k$ are nonreal complex conjugate pairs then $\sigma_j = \sigma_k.$
We denote the set of admissible vectors by $\mathbb{R}_{a}^{n}.$ It is an open cone and admits the
partial order $\sigma \leq \xi$ if and only if $\sigma_j \leq \xi_j,\ 1 \leq j \leq n.$
\begin{defi}
\label{quasilattice}
The quasilattice corresponding to $\sigma \in \mathbb{R}_{a}^n$ is
$$
  \mathfrak{L}(\sigma) =
  \{\, (V^T\ell)_1 \, : \, \ell \in \mathbb{Z}^n \hbox{ and } |(V^T\ell)_j| < \sigma_j \hbox{ for all } 2 \leq j \leq n \, \}.
$$
\end{defi}
\begin{lem}
\label{QL1}
If $\sigma, \xi \in \mathbb{R}_{a}^n$ then
$$0 \in \mathfrak{L}(\sigma) \hbox{ and } \mathfrak{L}(\sigma) = - \mathfrak{L}(\sigma),$$
$$\mathfrak{L}(\sigma) \subseteq \mathfrak{L}(\xi) \hbox{ if and only if } \sigma \leq \xi,$$
$$\mathfrak{L}(\sigma) + \mathfrak{L}(\xi) = \mathfrak{L}(\sigma + \xi).$$
\end{lem}
{\bf Proof} The first assertion is obvious and the second and third follow from the density
assertion of Lemma \ref{BSW}.
\begin{lem}
\label{QL2}
Let $\sigma \in \mathbb{R}_{a}^n.$ If $\ell \in \mathfrak{L}(\sigma) \backslash \{0\}$ then
$|(V^T\ell)_1| \geq \prod_{j=2}^n \sigma_{j}^{-1}.$
The minimal distance between the points in $\mathfrak{L}(\sigma)$ is
$2^{1-n}\prod_{j=2}^n \sigma_{j}^{-1}.$
\end{lem}
{\bf Proof} If $\ell \in \mathbb{Z}^n \backslash \{0\}$ then
$\prod_{j=1}^n (V^T\ell)_j \in \mathbb{Z} \backslash \{0\}$ is nonzero and a symmetric polynomial in
$\mathbb{Z}[\lambda_1,...,\lambda_n]$ so is a nonzero integer, implying the first assertion. The second assertion follows since lemma \ref{QL1} implies
that the differences between points in $\mathfrak{L}(\sigma)$ are in $\mathfrak{L}(2\sigma).$
\begin{lem}
\label{QL3}
If $\sigma \in \mathbb{R}_{a}^n$ and $L > |\det V| \, \prod_{j=2}^n \sigma_{j}^{-1}$ then $\mathfrak{L}(\sigma)$ contains a nonzero point in the interval $(-L,L).$
\end{lem}
{\bf Proof} Let $X = (V^T)^{-1}\left[(-L,L) \times (-\sigma_1,\sigma_1) \cdots (-\sigma_n,\sigma_n)\right]^T.$ Then $X = -X$ and
$$\hbox{volume } X = 2^n \, |\det V|^{-1}\, L\, \prod_{j=2}^n \sigma_j > 2^n.$$
The result follows from Minkowski's Lemma \ref{minkowski} since
if $\ell \in \mathbb{Z}^n$ then
$(V^T\ell)_{1} \in (-L,L) \cap \mathfrak{L}(\sigma)$ if and only if
$
    \ell \in X.
$
\begin{lem}
\label{QL4}
If $\sigma \in \mathbb{R}_{a}^n$ then there exists $M(\sigma) > 0$ such that every open interval of length $M(\sigma)$ contains a point in $\mathfrak{L}(\sigma).$
\end{lem}
{\bf Proof} Let $G = \mathbb{R}^n/(VV^T\mathbb{Z}^n).$ Since $VV^T$ is
nonsingular and has integer entries $VV^T \mathbb{Z}^n$ is a rank $n$
subgroup of $\mathbb{Z}^n$ so $G$ is isomorphic to a torus group. Let
$q_n : \mathbb{R}^n \rightarrow G$
and
$\varphi_n : G \rightarrow \mathbb{T}^n$
be the canonical epimorphisms and observe that
\begin{equation}
\label{equivariant}
    p_n\left((VV^T)^{-1}x\right) = \varphi_n \circ q_n(x), \ \ x \in \mathbb{R}^n.
\end{equation}
Since the entries of $v = [\, 1, \, \lambda, \, ... \, , \lambda^{n-1}\, ]^T$ are linearly independent over $\mathbb{Q}$ the density assertion in Lemma \ref{BSW} implies that
\begin{equation}
\label{density1}
    G = q_n\left( (V
    \left[\, \mathbb{R} \times (-\sigma_1,\sigma_1) \times (-\sigma_n,\sigma_n) \, \right]^T \right).
\end{equation}
Since $G$ is compact there exists $\kappa >0$ such that
\begin{equation}
\label{density2}
    G = q_n\left( (V
    \left[\, (-\kappa,\kappa) \times (-\sigma_1,\sigma_1) \times (-\sigma_n,\sigma_n) \, \right]^T \right).
\end{equation}
Therefore Equation \ref{equivariant} implies that
\begin{equation}
\label{density3}
    \mathbb{T}^n = \varphi_n(G) = p_n\left( (V^T)^{-1}
    \left[\, (-\kappa,\kappa) \times (-\sigma_1,\sigma_1) \times (-\sigma_n,\sigma_n) \, \right]^T \right)
\end{equation}
and the result follows by choosing $M(\sigma) = 2\kappa.$
\begin{lem}
\label{QL5}
If $|c_0| = 1$ then
$\lambda \mathfrak{L}(\sigma) = \mathfrak{L}([\, 0,\, |\lambda_2|\, \sigma_2,...,\, |\lambda_n|\, \sigma_n\, ]^T).$
\end{lem}
{\bf Proof}
Equation \ref{keyeq} gives $\lambda_j(V^T\ell)_j = (V^TC^T\ell)_j, \ j = 1,...,n$
and hence
$$
\lambda \mathfrak{L}(\sigma) =
    \{ \, (V^TC^T\ell)_1 \, : \, \ell \in \mathbb{Z}^n \hbox{ and }
        |(V^TC^T\ell)_j| < |\lambda_j|\, \sigma_j \hbox{ for } j = 2,...,n\, \}.
$$
Since $\det C = (-1)^n c_0 = \pm 1,$ $C^T \mathbb{Z}^n = \mathbb{Z}^n$ and hence
$\lambda \mathfrak{L}(\sigma) = $
$$
        \{ \, (V^T \, \ell)_1 \, : \, \ell \in \mathbb{Z}^n \hbox{ and }
        |(V^T\, \ell)_j| < |\lambda_j|\, \sigma_j \hbox{ for } j = 2,...,n\, \} =
        \mathfrak{L}([\, 0,\, |\lambda_2|\, \sigma_2,...\, ,|\lambda_n|\, \sigma_n\, ]^T).
$$
\begin{theo}
\label{multires}
If $|c_0| = 1$ and $\sigma \in \mathbb{R}_{a}^n$ then
$\xi = [\, 0,\, (1-|\lambda_2|)\, \sigma_2,...\, ,(1-|\lambda_n|)\, \sigma_n\, ]^T \in \mathbb{R}_{a}^n$
and if $f$ is a refinable distribution satisfying Equation \ref{refinable}
and $\tau_j \in \mathfrak{L}(\xi)$ then the spaces
\begin{equation}
\label{W}
    W_k = \hbox{span } \{\, f(\lambda^k\, x - \tau)\, : \, \tau \in \mathbb{L}(\sigma)\, \},
    \ \ k \in \mathbb{Z},
\end{equation}
satisfy $W_k \subset W_{k+1}.$
\end{theo}
{\bf Proof} If $\tau \in \mathbb{L}(\sigma)$ then Equation \ref{refinable} implies
that
$$f(\lambda^k \, x - \tau) = \sum_{j=1}^m a_j f(\lambda^{k+1} \, x - \lambda \tau - \tau_j).$$
Since $|c_0| = 1,$ Lemma \ref{QL5} implies that
$$
\lambda \, \tau \in \mathfrak{L}([\, 0,\, |\lambda_1|\, \sigma_1,...,\, |\lambda_n|\, \sigma_n\, ]^T),
$$
and hence Lemma \ref{QL1} implies that
$$
\lambda\, \tau + \tau_j \in
\mathfrak{L}([\, 0,\, |\lambda_1|\,\sigma+\xi_1,...\, ,|\lambda_n|\, \sigma_n+\xi_n\, ]^T) = \mathfrak{L}(\sigma).
$$
Theorem \ref{multires} provides the foundation of a multi-scale theory of refinable distributions
whose dilations are PV numbers $\lambda$ that are units in the ring of integers in the number field $\mathbb{Q}(\lambda)$ and whose translates are in quasilattices associated to $\lambda.$ The union
of all the quasilattices equals the set
$\mathbb{Z}[\, 1, \, \lambda, \, ... \, , \lambda^{n-1}\, ]$
which is a dense rank $n$ subset of $\mathbb{R}.$ It can be shown that the Fourier transform of the sum of unit point measures at each point of a quasilattice is a tempered distribution supported on the countable set
$\mathbb{Z}[\, 1, \, \lambda, \, ... \, , \lambda^{n-1}\, ].$
These quasilattices constitute a class of quasicrystals which have been extensively studied. We refer the reader to papers by Bombieri \cite{bombieri1}, \cite{bombieri2}, the book by Senechal \cite{senechal}, and references therein.

\Finishall
\end{document}